\documentclass[3p]{elsarticle}

\usepackage{amssymb}
\usepackage{amsmath}
\usepackage{amsthm}
\usepackage[all]{xy}
\usepackage{enumerate}
\usepackage[english]{babel}
\usepackage{multirow}
\usepackage{graphicx}
\usepackage{float}

\makeatletter

\newcommand{\Rmnum}[1]{\expandafter\@slowromancap\romannumeral #1@}
\makeatother

\begin{document}

\begin{frontmatter}
\author[rvt,focal]{Hanbaek Lyu}
\author[rvt,focal]{Piotr Grzegorz Jablonski}
\address[rvt,focal]{Seoul National University}

\title{\textbf{Four-Dimensional Discrete-time Lotka-Volterra Models with an Application to Ecology} }

\theoremstyle{definition}	\newtheorem{ex}{Example}[subsection]
\theoremstyle{definition}	\newtheorem{problem}{Problem}[subsection]
\theoremstyle{definition}	\newtheorem{proposition}{Proposition}[subsection]
\theoremstyle{definition}	\newtheorem{lemma}{Lemma}[subsection]
\theoremstyle{definition}	\newtheorem{corollary}{Corollary}[subsection]
\theoremstyle{definition}	\newtheorem{definition}{Definition}[subsection]
\theoremstyle{definition}	\newtheorem{theorem}{Theorem}[subsection]
\theoremstyle{definition}	\newtheorem{remark}{Remark}[subsection]
\theoremstyle{definition}	\newtheorem{fact}{Fact}[subsection]
\theoremstyle{definition}	\newtheorem*{introduction}{Introduction}
\theoremstyle{definition}	\newtheorem*{note}{Note}

\begin{abstract}
	This paper presents a study of the two-predators-two-preys discrete-time Lotka-Volterra model with self-inhibition terms for preys with direct applications to ecological problems. Parameters in the model are modified so that each of them has its own biological meaning, enabling more intuitive interpretation of biological conditions to the model . Moreover, the modified version is applicable to simulate a part of a large ecosystem, not only a closed predator-prey system with four species. An easy graphical method of analysis of the conditions on parameters ensuring long persistence under coexistence is presented. As an application, it is explained that a predator specie who feed on relatively small number of preys compared to the other predator species must be selective on the available preys in order for long persistence of the ecosystem. This may be regarded as a theoretical explanation of the existence of flush-pursuer birds, those who uses highly specialized hunting strategy and cross-adapts to the ecosystem relative to the ordinary bird species.
\end{abstract}

\end{frontmatter}

\section{Introduction}

This paper is a study of four-dimensional discrete-time Lotka-Volterra model by the method of linearization and simulations, with application to a realistic ecological system. The Lotka-Volterra models, proposed by Lotka \cite{Lotka} for the description of chemical reactions and by Volterra \cite{Volterra} for the description of interacting populations, have been the starting point in many dynamical models. The original and most commonly studied version of this model describes a continuous-time dynamics, whereas in real ecosystem the changes in populaions of each species due to competitive interaction cannot occur continuously. Hence it would be more suitable to use a discrete-time dynamical system, represented by difference equations or, more properly, by the iterated application of maps, are often a more suitable tool for modeling the dynamics in ecosystem. 

The choice of discrete-time dynamic models is crucial because they may exhibit more complicated dynamic behaviors then continuous-time models. Indeed, even in one-dimensional discrete dynamical systems represented by iterated quadratic maps, like the well known logistic map, periodic and chaotic trajectories can be easily observed \cite{Period3}. Similarly, two-dimensional discrete dynamical systems, in particular discrete time Lotka–Volterra equations, can exhibit a wide variety of complicated asymptotic behaviors, from convergence to a fixed point or a periodic cycle until quasi-periodic motion along a closed invariant curve or even an erratic motion inside a bounded chaotic attractor (see e.g. \cite{Blackmore}, \cite{Liucomplex}, \cite{Liubifurcation}).

The local asymptotic behaviour near a fixed point of the system can be analyzed by linearization, which is basically an approximation of the system by a linear system defined by the Jacobian at the given fixed point. The global behaviour of a discrete-time dynamical system, however, is far more interesting and complicated. One of a few methods available for the study of it is the elegant theory Lyaponov; if one can find a scalar function that satisfies certain conditions, called the Lyaponov function, then one can say about the global asymtotic behaviour to certain points or bounded in certain region. But since their is no systematic way to find the Lyapunov function of a given dynamical system, this is often a hard task. The method of critical sets is another way of studying the global behaviour of discrete-time dynamical system which is applicable if the transition map is non-invertible, as it happens to be in the case of Lotka–Volterra maps. Such method is a powerful tool introduced in the Seventies, by which the global dynamical properties, especially the bifurcations, can be usefully characterized(see \cite{Mira}, \cite{Gumowski}). See \cite{Bischi} for such a study of global properties of three dimensional Lotka-Volterra map by means of critical set. But the method of critical set applied to a map of dimension greater than two, may contain a very hard and challenging task, due to the difficulties met in the computer-assited graphical visulaization(see e.g. \cite{Bischi3D}, \cite{Hauser}) There is also an excelent mathematical analysis of a special kind of three-dimensional Lotka-Volterra map by Gardini et al. \cite{Gardini}

However, those methods are too difficult to be practical for biologists who wants to use the mathematical modeling to an actual biological problem in order to support an emperical hypothesis or make a theoretical hypothesis to guide further emperical study. To present a sound and practical modeling stretage which is easily applicable to some interesting biological problems, we need to satisfy the following two points. First, the mathematical model must reasonably simulate the biological systems in interest. Second, the method of analyzing such a mathematical system need to be handy and could draw a meaningful biological implications. In the present paper, we concern a four dimensional discrete-time Lotka-Volterra system with 2 preys - 2 predators
\begin{equation}
\begin{cases}
	x_{1}(t+1)=x_{1}(t)[r_{1}-k_{1}x_{1}(t)-a_{11}X_{1}(t)-a_{12}X_{2}(t)]\\
	x_{2}(t+1)=x_{2}(t)[r_{2}-k_{2}x_{2}(t)-a_{21}X_{1}(t)-a_{22}X_{2}(t)]\\
	X_{1}(t+1)=X_{1}(t)[-p_{1}+b_{11}x_{1}(t)+b_{12}x_{2}(t)]\\	
	X_{2}(t+1)=X_{2}(t)[-p_{2}+b_{21}x_{1}(t)+a_{22}x_{2}(t)]
\end{cases}
\label{dLVsystem}
\end{equation}
which we consider as the one satisfying the two requirements in some acceptable extent. First, the system in concern can simulate biological situations : 1) a system where a specie of generalist predator and a specie of specialist predator which is selective at a certain prey specie competes over two prey species; 2) a system with a predator specie heavily dependent on two prey species and with another predator specie which is less dependent one the two prey species due to the presence of other source of food; 3) a closed system with single prey and predator species to which a new predator specie is introduced by migration or mutation. 4) a closed system with single prey and predator species to which a new prey specie is introduced by migration or mutation. Second, as a method of analyzing such system, we introduce an application of linearization which enables us to draw some biologically meaningful conclusions about the global behaviour of the system and conditions on prey and predator species for coextence in the current system. 

\section{The two-dimensional Lotka-Volterra map}
\section{The three-dimensional Lotka-Volterra map}
The three-dimensional Lotka-Volterra map $T_{3}:(x_{1},x_{2},x_{3})\mapsto (x_{1}',x_{2}',x_{3}')$ is defined by 
\begin{equation}
	T_{3} : 
\begin{cases}
	x_{1}'=x_{1}(e_{1}+a_{11}x_{1}+a_{12}x_{2}+a_{13}x_{3})\\
	x_{2}'=x_{2}(e_{2}+a_{21}x_{1}+a_{22}x_{2}+a_{23}x_{3})\\
	x_{3}'=x_{3}(e_{3}+a_{31}x_{1}+a_{32}x_{2}+a_{33}x_{3})
\end{cases}
\end{equation}

\section{The four-dimensional Lotka-Volterra map}
The general four-dimensional Lotka-Volterra map $T_{4}$ can be defined as similarly as the three dimensional map $T_{3}$. However, we concentrate on the following specialization $T_{4}':(x_{1},x_{2},X_{1},X_{2})\mapsto (x_{1}',x_{2}',X_{1}',X_{2}')$ defined by 
\begin{equation}
	T_{4}' : 
\begin{cases}
	x_{1}'=x_{1}(r_{1}-k_{1}x_{1}-b_{11}X_{1}-b_{12}X_{2})\\
	x_{2}'=x_{2}(r_{2}-k_{2}x_{2}-b_{21}X_{1}-b_{22}X_{2})\\
	X_{1}'=X_{1}(-p_{1}+c_{11}x_{1}+c_{12}x_{2})\\	
	X_{2}'=X_{2}(-p_{2}+c_{21}x_{1}+c_{22}x_{2})
\end{cases}
\label{dVLmap}
\end{equation}
where all the parameters are positive real numbers. The discrete-time dynamical system determined by the repeated application of this map $T_{4}'$ simulates a predator-prey ecosystem consisting of two prey and predator species; the lower and upper case variables denote prey and predator species, repectively. The underlying biological assumptions for this model are as follows.

\begin{description}
	\item{(a).} All species have overlapping generations; the populations of all species are updated simultaneously.
	\item{(b).} There is no direct interaction within preys or predators. 	
	\item{(c).} We does not consider any age structure in the model. Any two individuals of the same species are equivalent(have the same traits) regardless of their age.
	\item{(d).} Without the presence of predators, the prey populations follow logistic growth.
	\item{(e).}  Without the presence of preys, the predator populations follow exponential decay.
	\item{(f).}  The amount of predation of prey $i$ by predator $j$ and reproduction of predator $j$ due to the consumption of prey $i$ at each generation are propotional to the product $x_{i}X_{j}$ of their population, as such product represents the number of encounters during that generation. 
\end{description}

Analyzing the map \eqref{dVLmap}, it is easy to see that the coordinate axes, planes, and spaces are invariant sets, i.e., if a coordinate of $\mathbf{x}\in \mathbb{R}^{4}$ is zero, then the same coordinate of $T_{4}'(\mathbf{x})$ is zero. For example, a trajectory starting from a point such that $x_{1}(0)=0$ remains trapped inside the 3-dimensional space $x_{1}=0$ and is determined by the restriction of \eqref{dVLmap} obtained from setting $x_{1}\equiv 0$. One easily realizes there are total 15 fixed points with at least one vanishing coordinate. Now suppose $\mathbf{x}_{*}$ is a \textit{non-degenerate} fixed point of \eqref{dVLmap}, i.e., a fixed point without zero coordinate. It is clear that such a fixed point is given as a solution of the four linear equations obtained form canceling out the degree 1 factors in each of the four equations in \eqref{dVLmap}. In other words, we have 
\begin{equation}
	\mathbf{x_{*}}=
		\begin{bmatrix} 
			r_{1} \\ r_{2} \\ -p_{1} \\ -p_{2} 
		\end{bmatrix}
		\begin{bmatrix} 
			k_{1} & 0 & b_{11} & b_{12} \\
			0 & k_{2} & b_{21} & b_{22} \\
			-c_{11} & -c_{12} & 0 & 0 \\	
			-c_{21} & -c_{22} &  0 & 0 
		\end{bmatrix}^{-1}
\label{fixedpoint}
\end{equation}
Of course, the $4\times 4$ matrix inverse in \eqref{fixedpoint} may not exist and there could be no or many non-degenerate fixed points of \eqref{dVLmap}. Such cases are not of our interest here, however, and as long as we are interested in modeling ecosystems, there is no serious problem in assuming that such inverse matrix exists since the probability that the inverse matrix does not exist when the parameters of \eqref{dVLmap} are arbitrarily given, is zero. Hence we may assume that \eqref{dVLmap} has a unique non-degenerate fixed point given by  \eqref{fixedpoint}. We will be more interested in such a non-degenerate fixed point if all of its four coordinates are positive. We may call such fixed point \textit{positive}. 

After all, we are interested in the condition for the parameters of the map $T_4$ such that it has a unique positive equilibrium and the system persists near the fixed point for a relatively long time. This is the situation when the corresponding predator-prey system persists with all the four species coexisting. For this task, we shall use the standard technique of \textit{linearization}, which we will briefly summarize. 

As one studies the trend of a curve near a certain point by looking at the tangent line, one can study a non-linear system near a fixed point by approximating to a well-understood relatively simple system, which is, so to say, “tangent” to the non-linear system at the given fixed point. Such a simple system is called linear system, whose transition map is given by matrix multiplication. For example, if A be a $4\times 4$ matrix with real entries, then we have a four-dimensional linear system given by the transition map 
\begin{equation}
	\mathbf{x}'=A\mathbf{x}. \label{linearmap}
\end{equation}
Thus gien an initial condition $\mathbf{x}(0)$, the state $\mathbf{x}(n)$ after $n$ iteration of the map \eqref{linearmap} is given by 
\begin{equation}
	\mathbf{x}(n)=A^{n} \mathbf{x}(0). \label{lineartrajectory}
\end{equation}
If, in particular, the initial condition $\mathbf{x}(0)$ is an eigenvector of $A$ with eigenvalue $\lambda$, then the trajectory will be 
\begin{equation}
	\mathbf{x}(n)=\lambda^{n} \mathbf{x}(0). \label{particularsolution}
\end{equation}
Hence $\mathbf{x}(n)$ converges to the origin if $|\lambda|<1$, and diverges from the origin if $|\lambda|>1$, in both cases with speed $|\lambda|$. If one requires all the trajectories converge to the origin, it means that the matrix $A^n$ converges to the zero matrix. In this case we say that A is zero-convergent and the origin is a stable fixed point of \eqref{linearmap}. If some trajectory tends to infinity, then the origin is said to be unstable. More loosely, if one requires that every trajectory is trapped in certain region, it means that there is a positive real number $M$ such that every entry of the matrix $A^n$ has absolute value less than $M$ for all $n\ge 1$. In this case we call $A$ power-bounded. Evidently, a matrix $A$ is zero-convergent only if every eigenvalue must has modulus less than 1, and power-bounded only if every eigenvalue has modulus no greater than 1. It is well-known that these conditions are also sufficient. This is the condition for long persistence of a discrete linear system. Hence one can determine the dynamics of a linear system by looking at the eigenvalues of the transition matrix.

On the other hand, one can find such a linear system which is “tangent“ to a given non-linear system at a given fixed point by finding the Jacobian matrix of the transition map at the fixed point. In our case, let $\mathbf{x}_{*}=(x_1^*,x_2^*,X_1^*,X_2^*)$ be a fixed point of the map $T_4'$. By using the equilibrium condition, the Jacobian matrix $J_4'$ of $T_4'$ at $\mathbf{x}_{*}$ is given by   
\begin{equation}
	J_{4}'=\begin{bmatrix}
	1-k_{1}x_{1}^{*} & 0 & -a_{11}x_{1}^{*} & -a_{12} x_{1}^{*} \\
	0 & 1-k_{2}x_{2}^{*} & -a_{21}x_{2}^{*} & -a_{22} x_{2}^{*} \\
	b_{11}X_{1}^{*} & b_{12}X_{1}^{*} & 1 & 0 \\
	b_{21}X_{2}^{*} & b_{22}X_{2}^{*} & 0 & 1 
\end{bmatrix}
\label{Jacobian}
\end{equation}
The derived linear system 
\begin{equation}
	\mathbf{x}'=J_{4}'\mathbf{x}. \label{jacobiansystem} 
\end{equation}
is the one to which we approximate the dynamics of our system near the fixed point $\mathbf{x}_{*}$. As we have seen, the system \eqref{jacobiansystem} is locally stable at the fixed point $\mathbf{x}_{*}$ if the four eigenvalues of $J_4'$ is less than 1 in modulus, and unstable if at least one eigenvalue has modulus greater than 1. Moreover, one can guess that in order for the trajectories near the fixed point to be bounded, the Jacobian $J_4'$ must be at least power-bounded. In other words, if the system \eqref{jacobiansystem} shows a long persistence near a fixed point, then it is very much likely that the eigenvalues of the corresponding Jacobian matrix is no greater than 1 in modulus. Even though there are not easy algebraic manipulations to express the four eigenvalues of $J_4'$ in terms of the parameters, the formulae  (4) and (5) enables us to determine whether the system \eqref{jacobiansystem} with a particular set of parameters has a stable non-degenerate fixed point or not. 

\subsection{Modeling}

Although the map $T_4'$ is handy for mathematical analysis, it is not so practical to use to simulate a particular type of ecosystem; it is not clear how to assign suitable values of the parameters to reflect important characteristics of the subject ecosystem. Look at the equations in \eqref{dVLmap}. The parameters $r$ and $k$ describes the intrinsic population dynamics of preys, and $p$ describes the dependency of predators to preys. What describes the interaction between the four species are the parameters b and c, denoting the predation rate of preys and reproduction rate of predators, respectively. However, these two parameters are the results of combined effect of various factors between predators and preys. We assume the predation of preys is determined by the following three factors:
\begin{description}
	\item{(1)}	\textit{hunting efficiency} of a particular predator specie over a particular prey specie;
	\item{(2)}	\textit{adaptation} of a particular prey specie to a particular predator specie;
	\item{(3)}	\textit{search rate} of a particular predator specie.
\end{description}
Let $e_ji$ be the hunting efficiency of predator $j$ over prey $i$, $a_ji$ be the adaptation coefficient of prey $i$ to predator $j$, and $s_i$ be the search rate of predator $i$. We assume the number of encounters between predator $j$ and prey $i$ is $s_i x_j X_i$, and the number of prey $j$ predated by predator $i$ during one generation is given by 
\begin{equation}
	s_i a_{ij} e_ji x_i X_j
\end{equation}

Note here that the larger the value of $a_{ij}$, the more the prey $i$ is adapted against predator $j$. Hence we may factor the predation rate $b_{ij}$ by $b_{ij}=s_i a_{ij} e_{ji}$. On the other hand, the amount (10) of hunted prey $i$ will be consumed and contribute to the reproduction of predator $j$ for the next generation. We multiply another coefficient $q_{ji}$, the conversion ratio, to (10), to measure the amount of reproduction of predator $j$ due to the consumption of prey $i$. Hence we may factor the reproduction rate $c_{ji}$ by $c_{ji}=q_{ji} a_{ij} e_{ji}$. Finally, we may write $k_i=r_i/K_i$, where $K_i$ is the carrying capacity of prey $i$. Now we reformulate the map $T_4'$ as follows: 
\begin{equation}
	T_{4}' : 
\begin{cases}
	x_{1}'=r_{1}x_{1}(1-\frac{x_{1}}{K_{1}})-s_{1}a_{11}e_{11}x_{1}X_{1}-s_{2}a_{12}e_{21}x_{1}X_{2}\\
	x_{2}'=r_{2}x_{1}(1-\frac{x_{2}}{K_{2}})-s_{1}a_{21}e_{12}x_{2}X_{1}-s_{2}a_{22}e_{22}x_{2}X_{2}\\
	X_{1}'=X_{1}-p_{1}X_{1}+s_{1}q_{11}a_{11}e_{11}x_{1}X_{1}+s_{1}q_{12}a_{21}e_{12}x_{2}X_{1}\\
	X_{2}'=X_{2}-p_{2}X_{2}+s_{2}q_{21}a_{12}e_{21}x_{1}X_{2}+s_{2}q_{22}a_{22}e_{22}x_{2}X_{2}
\end{cases}
\label{dLVmodeling}
\end{equation}
This model will not only be interpreted as a model of a closed predator-prey system with four species, but also a model of a portion of a large ecosystem, according to the values of the parameters and their biological meaning.  

Figure 1 illustrates six numerical simulations that show various dynamical patterns of the system \eqref{dVLmodeling}. The population of the four species are plotted versus the number of generations elapsed. Red, green, black, blue dots refer to $x_1,x_2,X_1,X_2$, respectively. 
\begin{figure}[H]
	\begin{center}
		\includegraphics[width=15cm,height=10cm]{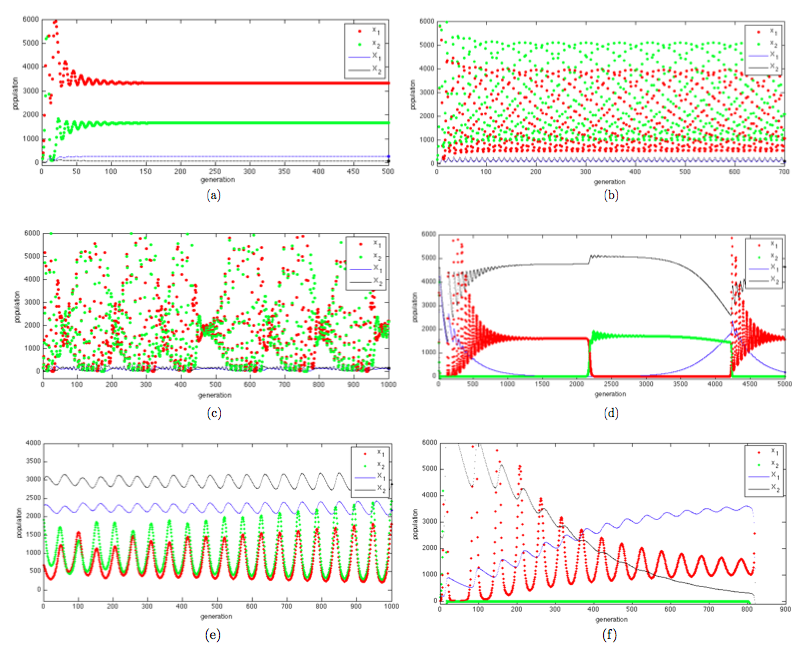}
	\end{center}
	\caption{(a), (b), (c) have common parameters $r_{1}=r_{2}=1.5$, $K_{1}=K_{2}=10^{4}$, $s_{1}=s_{2}=0.01$, $p_{1}=p_{2}=0.3$, and $q_{11}=q_{12}=q_{21}=q_{22}=0.02$ but different hunting rate; (a) $a_{11}=a_{12}=0.3$, $a_{21}=0.2$, $a_{22}=0.5$  (b) $a_{11}=0.2$, $a_{12}=0.45$, $a_{21}=0.5$, $a_{22}=0.25$  (c) $a_{11}=0.24$, $a_{12}=0.6$, $a_{21}=0.56$, $a_{22}=0.27$. (d) has the same value of $r_{i}$s and $q_{ij}$s, but $K_{1}=K_{2}=10^{5}$, $p_{1}=p_{2}=0.01$, $s_{1}=s_{2}=0.001$, $a_{11}=0.2$, $a_{12}=0.4$, $a_{21}=0.31$, $a_{22}=0.29$.  (e) and (f) have the same parameters; their set of parameters differ from that of (d) only in the $a$'s; $a_{11}=0.4$, $a_{12}=a_{21}=0.2$, $a_{22}=0.35$. The first four have the same initial condition $(100,100,100,100)$, whereas $(700,2000,2300,3000)$ for $(e)$ and $(10,100,200,2500)$ for $(f)$. The non-degenerate fixed points are (a) $(3.3333e+3, 1.6667e+3, 277.7778, 83.3333)$, (b) $(1.7143e+3, 2.5714e+3, 140.8163, 192.2449)$, (c) $(1.8252e+3, 1.7699e+3, 132.8351, 162.0379)$, (d) $(833.3333, 833.3333, 450.7576, 4.5076e+3)$, (e), (f) (750.0000m 1.0000e+3, 2.2406e+3, 2.9625e+3)}  . 
\end{figure}
The first four examples starts at the same initial point, and the first three cases show different dynamics depending on the values of hunting rates. Hence the two prey species are identical with the same intrinsic rate of increase and carrying capacities, and the two predator species have the same dependency on preys and conversion ratios. These parameters are suitably chosen to allow long persistence. However, the strategies of predators are what make the difference. In (a), predator 1 is equally efficient at hunting both prey 1 and 2 with $a_{11}=a_{12}=0.3$, but predator 2 is specialised in hunting prey 2, with $a_{21}=0.2$ and $a_{22}=0.5$. Hence it simulates a system with a generalist and a specialist predator. The system converges to the unique positive fixed point, that is, there is a stable coexistence at the fixed population. 

In (b), both predators are specialized in different prey species. The dynamics is a regular oscillation around the fixed population, and the persistence is very strong. Example (c) also simulates a similar system, but now the predator 1 is slightly more specialized than before. The dynamics is a irregular oscillation around the fixed population, with characteristic pattern that the period of oscillation continuously increases and decreases alternately.

In (d) we have 10 times more carrying capacities for the preys and less dependencies of predators and search rates. Predator 1 is specialized in hunting prey 2, while predator 2 is almost equally efficient on both preys. Its dynamics is interesting, showing a domination by the generalist predator 2, with alternating dominance of prey 1 and 2. Under the dominance of predator 2, the number of predator 1 keeps decreasing and when this number goes down below a certain level, the prey 2 suddenly bursts out and prey 1 decreases rapidly. Due to the abundance of prey 2 to which the predator 1 is specialized, the number of predator 1 gradually increases. But at some point large number of predator 1 almost extinguishes the prey 2, and then prey 1 now thrives with rapid oscillation. 

The last two examples show the dependence on the initial conditions. Both (e) and (f) have the same parameters but different initial populations. The initial population for (e) is more close to the fixed population, and it shows almost regular oscillation around the fixed populations. On the other hand, the system of (f) undergoes increasing or decreasing oscillations and then suddenly collapses after about 600th generations. 

Nonetheless, all the six examples persist for a very large number of generations, mostly around the fixed populations. Such a long persistence around the fixed population would be not so much likely if there is an eigenvector with eigenvalue whose absolute value far exceeds 1, since then any component of the initial population vector along the eigenvector will be rapidly repelled from the fixed point. Indeed, the table below shows that the eigenvalues in the six examples are mostly less than 1, and some of them are slightly greater than 1.

Notice in the three unstable cases, the absolute values of the four eigenvalues are close to the critical value 1, especially the one exceeding 1. It seems that a long persistence of the system around the unstable fixed point is only possible if the \textit{spectral radius} of the Jacobian, i.e., the maximum value of the modulus of the eigenvalues, which exceeds 1, is closed to 1. The decaying oscillation toward the fixed point shown in (a) is due to the complex eigenvalues less than 1 in modulus. The presence of complex eigenvalues greater than 1 in modulus as in the other three unstable cases may contribute to the persisting oscillatory behaviors. However, it would not be possible to predict the global pattern of the system solely from the local data of eigenvalues. 
\begin{figure}[H]
	\begin{center}
		\includegraphics[width=12cm,height=4.7cm]{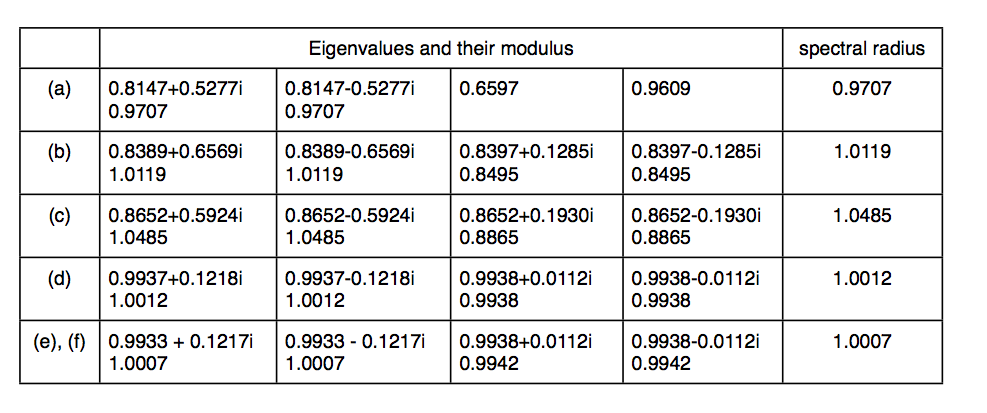}
	\end{center}
	\caption{Eigenvalues and their modulus of the corresponding Jacobian matrix}
\end{figure}

local repelling force, attracting force

global attracting force : self-inhibition term for preys inhibit total collapse

\subsection{Stability diagram analysis}

In this section, we introduce the notion of $\textit{stability diagram}$ for the system \eqref{dLVmodeling}. We may fix the parameters other then those hunting efficiencies $a_{ij}$, and investigate the quadraples $(a_{11},a_{12},a_{21},a_{22})$ for which the system \eqref{dLVmodeling} have stable positive equilibirum. But the result would be plotted in 4-dimensional space, which is impossible to visulaize. To reduce the dimensionality, we may simply assume that 
\begin{equation}
a_{11}+a_{12}=1,\quad a_{21}+a_{22}=1. \label{constraint}
\end{equation}
The biological assumption under the constraint \eqref{constraint} is that each predator species has a limited hunting ability to be distributed over two prey species. For example, if predator 1 can hunt down prey 1 at rate $0.7$, then the hunting rate for prey 2 must be $0.3$, the remainder of the former. But of course two predator species can have different total hunting ability, and this point is reflected in the \textit{search rate}. Namely, if predator 1 have search rate twice of that of predator 2, meaning that predator 1 encounters preys twice more than predator 2, then predator 1 actually have twice total hunting ability than that of predator 2. 

Moreover, the constraint \eqref{constraint} does not harm the generality of our model, since one can \textit{normalize} the parameters of a given system (8) by replaing each paratmeter $a_{ij}$ by $a_{ij}'$ and $s_{i}$ by $s_{i}'$, where 
\begin{equation*}
	a_{11}'=\frac{a_{11}}{a_{11}+a_{12}}, \quad a_{21}'=\frac{a_{21}}{a_{11}+a_{12}}, \quad a_{21}'=\frac{a_{21}}{a_{21}+a_{22}}, \quad a_{22}'=\frac{a_{22}}{a_{21}+a_{22}}
\end{equation*}
\begin{equation}
	s_{1}'=s_{1}(a_{11}+a_{12}),\quad  s_{2}'=s_{2}(a_{21}+a_{22}). 
\label{normalization}
\end{equation}
Notice that the "total hunting ability" of prey $i$ is reflected in the factor $(a_{i1}+a_{i2})$ of the normalized serch rate $s_{i}'=s_{i}(a_{i1}+a_{i2})$. Also, now it make sense to call the hunting rate as "hunting probability" since the normalized hunting rate $a_{ij}'$ has value $0\le a_{ij}\le 1$. We may call a given four-dimensional Lotka-Volterra map \eqref{dLVmodeling} \textit{normalized} if it satisfies the constraint \eqref{constraint}. 

By normalization, the tuple $(a_{11},a_{12},a_{21},a_{22})$ becomes $(a_{11},1-a_{11},1-a_{22},a_{22})$, and hence determined by the two parameters $a_{11}$ and $a_{22}$. Thus we can plot all the pairs $(a_{11},a_{22})$ in the unit squre $S=\{(x,y)\in \mathbb{R}^{2}\,|\, 0\le x, y\le 1\}$ for which the given normalized system has \textit{positive and stable fixed point}. We call the set of all such pairs $(a_{11},a_{22})$ \textit{the stable region}, and such diagram as the \textit{stability diagram} of the given normalized system. As we have observed in the simulations in section 4.1, a long persistence of a given system around the non-degenerate fixed point might be possible only if the spectral radius of the system's Jabobian is very closed to 1. Even though we are dealing with a discrete dynamical system, the spectral radius of Jacobian for a given system is a continuous function in the parameters. In other words, if the two tuples $(a_{11},a_{22})$ and $(a_{11}',a_{22}')$ in the stability diagram are very close, then the corrdsponding spectral radii will be closed as well. Hence, we obtain the following principle for long persistence: \textit{a system corresponding to the tuple $(a_{11},a_{22})$ would be persistant for a long periods of generations only if the tuple $(a_{11},a_{22})$ is closed to the stable region}.

Given the values of parameters of system \eqref{dLVmodeling}, one can readily determine whether the corresponding system has a positive stable fixed point algorithmically, using the formulae \eqref{fixedpoint} and $\eqref{Jacobian}$. In stepwise discription, 
\begin{description}
	\item{(a)} Find the nondegenerate fixed point $\mathbf{x}_{*}$ of the given system by using the formula \eqref{fixedpoint} and check if all the coordinates are positive. 
	\item{(b)} Using \eqref{Jacobian}, find the Jacobian matrix of the system at the fixed point $\mathbf{x}_{*}$ and calculate the four eigenvalues.   
	\item{(c)} $(a_{11},a_{22})$ is a stable tuple if all of the the four coordinates of fixed point $\mathbf{x}_{*}$ are positive, and the maximum of the four eigenvalues in step (b) have minimum absolute value less then 1. 
\end{description}
Therefore, one can easily write a program in Matlab that draws the stability diagram for given parameters except the hunting rates, by implementing the stepwise algorithm as described above. So this could be a practical and useful tool that visualizes the results of localization analysis, to "see" what combinations of $a_{11}$ and $a_{22}$ may lead to a long persistance of the system with respect to the other given parameters. We will use this method for various set of parameters that represent certain type of ecosystems observable in nature.  

\textbf{Localist vs. Globalist.} First, consider a predator-prey system where "globalist" predator and "localist" predator compete. Globalist predators feed on various source of food since they search for preys on large area, whereas localist are localized so that they have very limited source of food compared to globalist. We are going to emulate this situation via our 4-dimensional model. Suppose the localist feed on two prey species. Then they ought have high search rate and dependencies on the two prey species, whereas the globalist have low values of them. 
\begin{figure}[H]
	\centering
		\includegraphics[width=6cm,height=3cm]{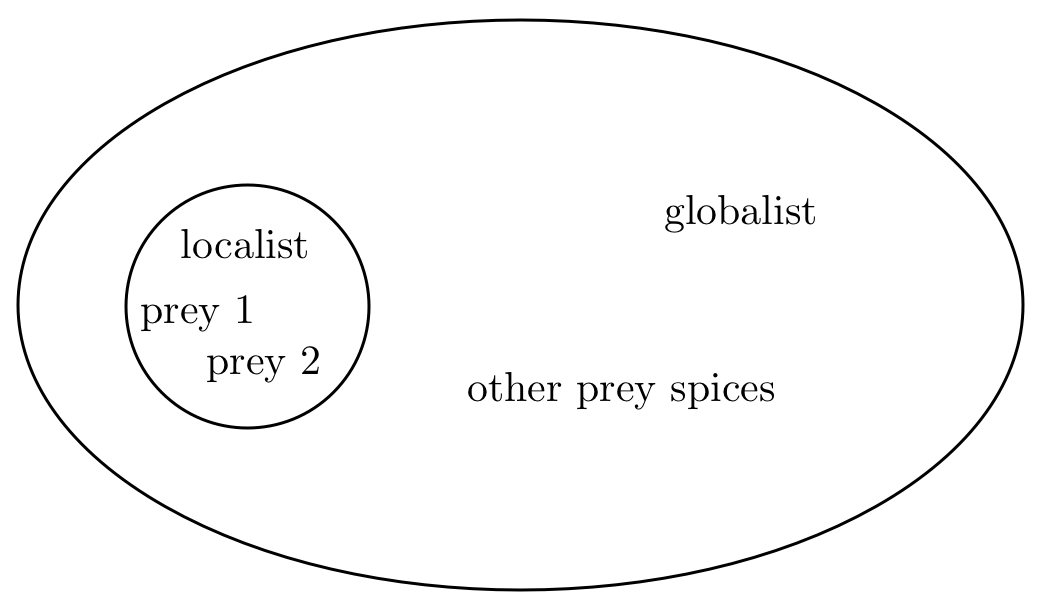}
	\caption{Habitat description of Example 5.0.1}
\end{figure}

Hence if we simulate the localist and globalist as predator species 1 and 2 in our model, then we would assume $s_{1}>s_{2}$ and $p_{1}>p_{2}$. We may at first leave the other parameters the same for the two predator species. In this situation, the stability diagram is as in the Figure 3.

\begin{figure}[H]
	\centering
		\includegraphics[width=13cm,height=9cm]{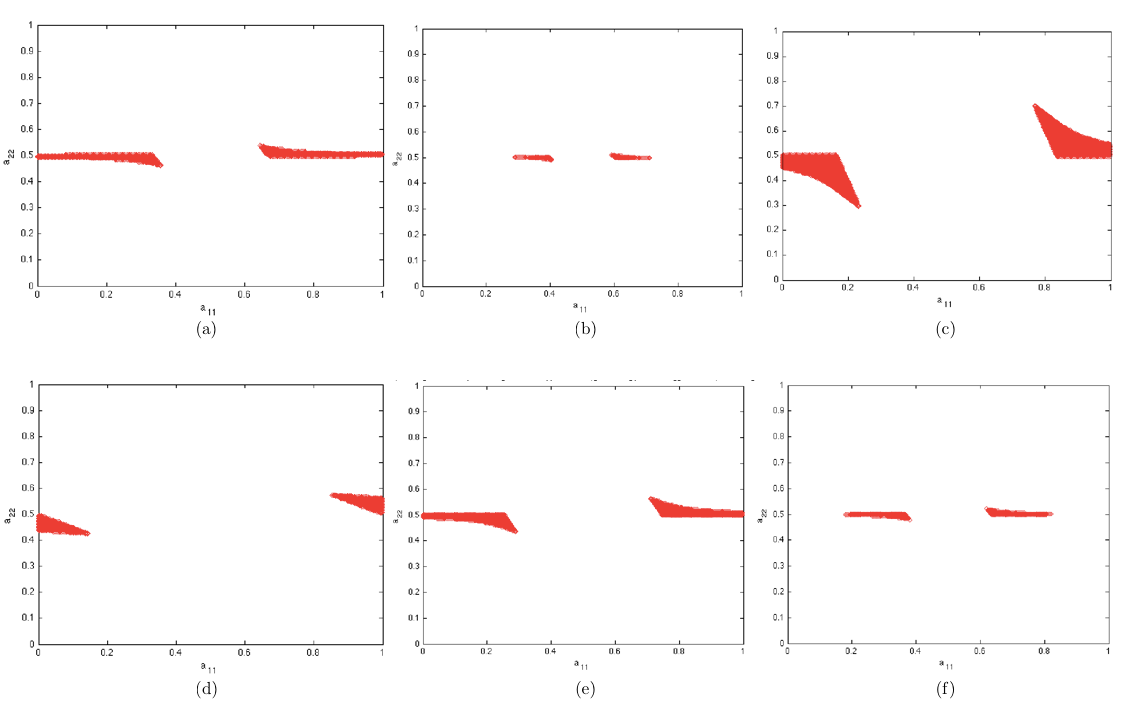}
	\caption{Stability diagrams for localist vs. globalist. The six cases have the common values of intrinsic rate of increse $r_{1}=r_{2}=1.5$ and carrying capacities $K_{1}=K_{2}=10^{5}$, and the first four has the same conversion ratios $q_{11}=q_{12}=q_{21}=q_{22}=0.01$. The other parameters for the first four are (a) $p_{1}=0.06$, $p_{2}=0.01$, $s_{1}=0.09$, $s_{2}=0.01$, (b) $p_{1}=0.06$, $p_{2}=0.01$, $s_{1}=0.09$, $s_{2}=0.012$, (c) $p_{1}=0.03$, $p_{2}=0.01$, $s_{1}=0.09$, $s_{2}=0.01$, (d) $p_{1}=0.06$, $p_{2}=0.01$, $s_{1}=0.03$, $s_{2}=0.01$. (e) and (f) have the same value of $p_{1}=0.06,p_{2}=0.01$, $s_{1}=0.09,s_{2}=0.01$ but they differ on the conversion ratios : (e) $q_{11}=a_{21}=0.13$, $q_{12}=q_{22}=0.01$,  (f) $q_{11}=q_{21}=0.01$, $q_{12}=q_{22}=0.011$.}
\end{figure}

Recall the red dots plotted on the diagram shows the tuples $(a_{11},a_{22})$ for which the system has positive and stable equilibrium. Hence if the system goes on stable state of coexistence, the tuple $(a_{11},a_{22})$ must be close to the red region. As one can see, the range of $a_{22}$ is very limited around 0.5, which means, the globalist must be generalist, that is, it should not be selective on the two local preys. On the other hand, the values of $a_{11}$ are apart from the center $0.5$. Hence the localist must be specialist, which means, selective on the two preys. Also consider the case (b) that generalist have slightly higher search rate $s_{2}=0.012$ from $0.01$ of (a). As one can see, in (b) the red regions ranges for $a_{11}$ have dramatically reduced down around the intervals $(0.3,0.4)$ and $(0.6,0.7)$. This means in this case the localist is not allowed to be too selective on preys. 

(c) is the stability diagram for lower dependency $s_{1}=0.03$ of localist. Now the localist must be more selective on preys, but the globalist can be little selective on preys. On the other hand, when the localist have lower search rate relative to its dependency, we have diagram (d), which shows pretty much restricted area of stable tuples $(a_{11},a_{22})$. Thus in any cases, we conclude that the localist are forced to be specialist, and globalist to generalist, in order for long coexistence of all the four species.

Now we consider other factors as well. First, suppose the localists are smaller than the globalists in body size. This difference will reflect in higher conversion ratios $q_{11}=q_{21}=0.013$ then the other $q_{21}=q_{22}=0.01$. But one might guess that the diagram in this case would be similar to that for Figure 14, since increasing both $q_{11}$ and $q_{21}$ woud have the same effect of increasing the search rate $s_{1}$. See (e). Similary, the case of large localist described by increased value of $q_{12}=q_{22}=0.011$, (f) which will have the same effect of increasing the search rate $s_{2}$, will likely to show a similar diagram as Figure 15. Hence regarding the body sizes, the interesting case will be asymmetric conversion ratios, say $q_{11}=0.01$, $q_{12}=0.011$, $q_{21}=0.011$ and $q_{22}=0.01$. 

\textbf{Adapted preys.} Assume that prey 1 is adapted to predator 1 and prey 2 to predator 2. We represent this situation by setting the adaptation coefficients $a_{11}=a_{22}=0.4$, $a_{12}= a_{21}=1.5$. We give four examples corresponding to different set of search rates s and predator dependency p.

\begin{figure}[H]
	\centering
		\includegraphics[width=12.5cm,height=12cm]{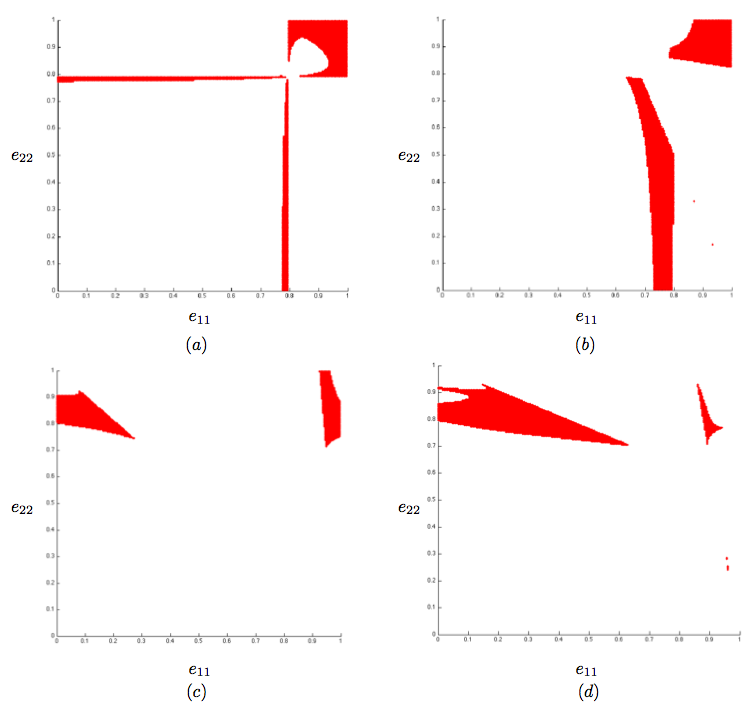}
	\caption{Common parameters : $r_1=r_2=1.5$, $K_1=K_2=10^5$, $q_{11}=q_{12}=q_{21}=q_{22}=0.01$, $a_{11}=a_{22}=0.4$, $a_{12}=a_{21}=1.5$. (a) $s_1=s_2=p_1=p_2=0.02$ (b) $s_1=s_2=p_2=0.02$, $p_1=0.03$ (c) $s_1=0.07$, $s_2=p_1=p_2=0.02$ (d) $s_1=0.07$, $p_1=0.04$, $s_2=p_2=0.02$.}
\end{figure}

In (a), all the parameters other than the adaptation coefficients are the symmetric in species 1 and 2. Hence the diagram is symmetric about the diagonal $y=x$. Ignoring the narrow strips near the line $e_{11}=0.8$ and $e_{22}=0.8$, the dominating upper right region shows that in order for long persistence of the system, predator $i$ should be far more efficient to prey $i$, which is much hard to catch than the other. Now we increase the dependency of predator 1 by 0.01 in (b). Now the vertical narrow strip is thicker, and the horizontal strip has been disappeared. As $e_{11}<0.8$, predator 2 has wide variety of options for its hunting strategy but not allowed to be too effective on hunting prey 2. On the other hand, as $e_{11}>0.8$, which means predator 1 is pretty much specialized for prey 1, then predator 2 must also be specialized for prey 2 in similar extent.  

Instead, if we increase the search rate $s_1$ from 0.02 to 0.07 in example (a), we get example (c). The upper-right region is the same, but we have upper-left region instead of the strap near the line $e_{11}=0.8$ in example (b). Hence a long persistence will be also possible if both the predator 1 and 2 are specialized for prey 2. 

Lastly, example (d) shows that increasing the value of $p_1$ in example (c) enlarges the upper left region, allowing wide range of options for predator 1. Comparing this example with (b), one sees that even if the predator 1 is more dependent on species, the coexistence condition might be more strict to the predator 2 if predator 1 has high search rate. 

\textbf{Small predators must be selective.}  Consider an ecosystem where two predators with different body size and two preys are competing. We want to investigate the stable tuples $(a_{11},a_{22})$ for which the system can persist with all the four species coexisting. To see this, set $q_{11}=q_{21}=0.04$ and $q_{21}=q{22}=0.01$, reflecting the predator specie 1 has small body size so that it has higher conversion ratio. Set other parameters symmetric for species 1 and 2. The stability diagram for this case is dipicted in Figure 5-(a), which looks like Figure 4-(a). As one can see, the stable values for $a_{11}$ are apart from the neutral 0.5; hence the small predator, specie 1, must be selective on preys. On the other hand, the other predator can be neutral on preys, but they are forced to prefer the other preys than those the small predator prefers. In fact, as we pointed out in previous examples, high $q_{11},q_{21}$ values does the same for high $s_{1}$ values. Hence small predators would be similar to the one who searchs preys more.

\textbf{More dependent on preys, more selective on preys.} Set $p_{1}=0.05$ and $p_{2}=0.01$, so that predator 1 is five times more dependent on the preys. The stability diagram Figure 5-(b) shows that Predator 1 must be selective on preys to in order for coexsistence.

\textbf{Assymetric diagrams.} The stability diagrams we have observed so far were all symetric about the point $(0.5,0.5)$, which means that there is no difference between the roll of species 1 and 2. Now let us consider more complex combination of parameters which yields assymetric stability diagrams. The diagram Figure 5-(c) shows predator 2 must be either specialized for prey specie 2 or neutral, while predator 1 can have any value of $a_{11}$, depending on $a_{22}$.The key factor that makes the diagram asymmetric in above example is the asymmetric values of conversion ratios $q_{ij}$. Such combinations of conversion ratios indicates that there are some non-trivial interrelationship between the four species other than just the body sizes. If the conversion ratio were dependent only on body sizes, for instance, we can assume prey 1 have body size 1 unit, so that predator 1 has size 100 unit since $q_{11}=0.01$. Similarly we would have that predator 2 has size $50=1/q_{21}$ unit, and hence prey 2 has size unit $3=100\cdot a_{21}$. However, the condition $q_{22}=0.01$ indicates that predator 2 has hundred times of body size of that of prey 2, which is not the case here. This is because the system of equations $z_{1}=az_{2}$, $z_{2}=bz_{3}$, $z_{3}=cz_{4}$ and  $z_{4}=dz_{1}$ may not have a solution in general. 

Figure 5-(d) is another example of such assymetric conversion ratios, with other parameters varied. This example describes a system where predator 1 searches more than predator 2,  while the former are more dependent on preys, and prey 2 are more reproductive, and conversion ratios are assymetric. As depicted in the diagram, the stable tuples are pretty much restricted than early examples. 

Lastly, Figure 5-(e) shows that such combination of parameters actually determines the stretagies of predators. 

\begin{figure}[H]
	\centering
		\includegraphics[width=13cm,height=9cm]{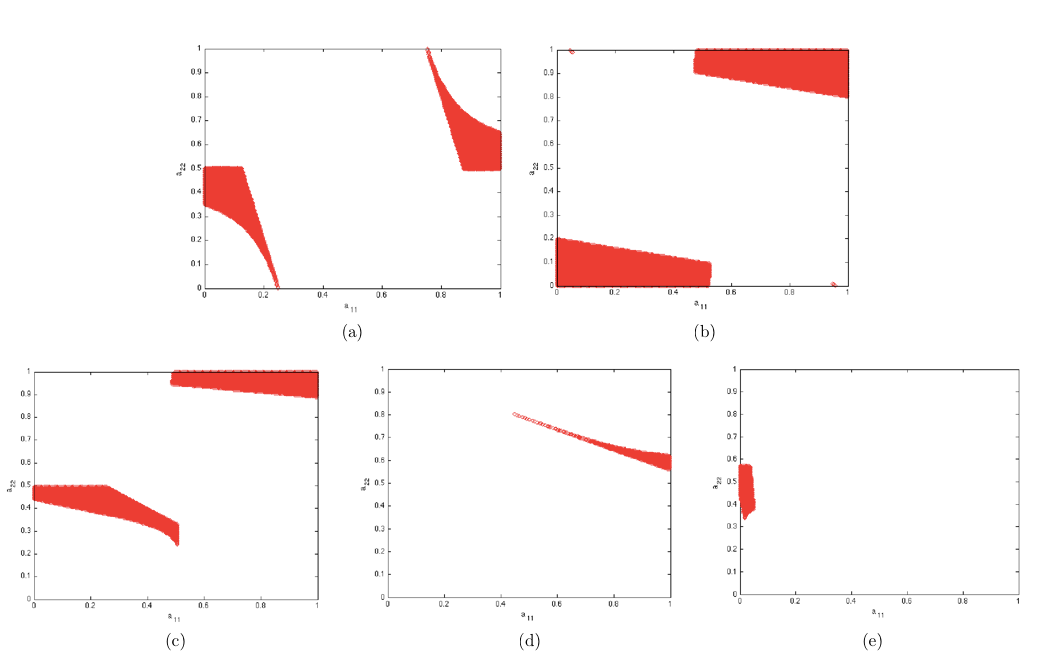}
	\caption{All the five cases have the same carrying capacities $K_{1}=K_{2}=10^{5}$, and the first three has the same value of $r_{1}=r_{2}=1.5$ and $s_{1}=s_{2}=0.01$. The other parameters are given as (a) $p_{1}=p_{2}=0.01$, $q_{11}=q_{21}=0.04$, $q_{12}=q_{22}=0.01$, (b) $p_{1}=0.05,p_{2}=0.01$, $q_{11}=q_{12}=q_{21}=q_{22}=0.01$, (c) $p_{1}=0.03,p_{2}=0.01$, $q_{11}=0.01,q_{12}=0.03,a_{21}=0.02,q_{22}=0.01$. (d) has $r_{1}=1.2,r_{2}=1.5$, $p_{1}=0.03$, $p_{2}=0.01$, $q_{11}=0.01,q_{12}=0.015,q_{21}=0.02,q_{22}=0.01$, $s_{1}=0.02$, $s_{2}=0.01$ and (e) has $r_{1}=1.3,r_{2}=1.7$, $p_{1}=p_{2}=0.01$, $q_{11}=0.03,q_{12=0.01},q_{21}=0.015, q_{22}=0.02$, $s_{1}=0.04, s_{2}=0.01$.}
\end{figure}

\section{Concluding remarks}

\section{Appendix}
\subsection{comparison of the four-dimensional discrete-time Lotka-Volterra system to the corresponding continuous-time system}
\subsection{comparison of the four-dimensional discrete-time Lotka-Volterra system to the discrete-system without self-inhibition term for prey species}
\subsection{definitions of oscilators and basin of oscillation}

\end{document}